\documentclass[11pt]{article}
\usepackage{amsfonts}
\usepackage{mathrsfs}
\usepackage{amsmath}
\usepackage{amssymb}
\usepackage{graphicx}
\usepackage{epsfig}
\usepackage{epic}
\usepackage{cite}
\usepackage{setspace}
\usepackage{amsthm,latexsym}
\usepackage{float}

\renewcommand{\paragraph}{\roman{paragraph}}
\newtheorem{theorem}{\scshape \mdseries \bf Theorem}[section]
\newtheorem{lemma}[theorem]{\scshape \mdseries  \bf Lemma}
\newtheorem{remark}[theorem]{\scshape \mdseries  \bf Remark}

\begin{document}

\title{\sf Full characterization of graphs having certain normalized Laplacian eigenvalue of  multiplicity $n-3$ }
\author{ \ Fenglei Tian\thanks{Corresponding author. E-mail address: tflqsd@qfnu.edu.cn.  Supported by '' the Natural Science Foundation of Shandong Province (No. ZR2019BA016)
''. }\ ,\ \
Yiju Wang \\
~~\\
\noindent{\small\it \ School of Management, Qufu Normal University, Rizhao, China.}
}
\date{}
\maketitle
\noindent {\bf Abstract:} \ Let $G$ be a connected simple graph of order $n$. Let $\rho_1(G)\geq \rho_2(G)\geq \cdots \geq \rho_{n-1}(G)> \rho_n(G)=0$ be the eigenvalues of the normalized Laplacian matrix $\mathcal{L}(G)$ of $G$.  Denote by $m(\rho_i)$ the multiplicity of the normalized Laplacian eigenvalue $\rho_i$. Let $\nu(G)$ be the independence number of $G$. In this paper, we give a full characterization of graphs with some normalized Laplacian eigenvalue of multiplicity $n-3$, which answers a remaining problem in [S. Sun, K.C. Das, On the multiplicities of normalized Laplacian eigenvalues of graphs, Linear Algebra Appl. 609 (2021) 365-385], $i.e.,$ there is no graph with $m(\rho_1)=n-3$ ($n\geq 6$) and $\nu(G)=2$. Moreover, we confirm that all the graphs with $m(\rho_1)=n-3$ are determined by their normalized Laplacian spectra.


\section{Introduction}

\quad { Throughout, only connected and simple graphs are discussed. Let $G=(V(G),E(G))$ be a graph with vertex set $V(G)$ and edge set $E(G)$. Let $N_G(u)$ be the set of all the neighbors of the vertex $u$. Then $d_u=|N_G(u)|$ is called the degree of $u$. For a subset $S\subset V(G)$, $S$ is called a set of twin points if $N_G(u)=N_G(v)$ for any $u, v\in S$. By $u\thicksim v$, we mean that $u$ and $v$ are adjacent. The distance of two vertices $u, v$ is denoted by $d(u,v)$ and the diameter of a graph $G$ is written as $diam(G)$.  A subset $S$ of $V(G)$ is called an independent set of $G$, if the vertices of $S$ induce an empty subgraph. The cardinality of the maximum independent set of $G$ is called the independence number, denoted by $\nu(G)$. The rank of a matrix $M$ is written as $r(M)$. Let $R_{v_i}$ be the row of $M$ indexed by the vertex $v_i$. The multiplicity of an eigenvalue $\lambda$ of $M$ is denoted by $m(\lambda)$.
Denote by $\mathcal{G}(n, n-3)$ the set of all $n$-vertex ($n\geq 5$) connected graphs with some normalized Laplacian eigenvalue of multiplicity $n-3$.
Let $A(G)$ and $L(G)=D(G)-A(G)$ be the adjacency matrix and the Laplacian matrix of graph $G$, respectively. Then the normalized Laplacian matrix $\mathcal{L}(G)=[l_{uv}]$ of graph $G$ is defined as $$\mathcal{L}(G)=D^{-1/2}(G)L(G)D^{-1/2}(G)=I-D^{-1/2}(G)A(G)D^{-1/2}(G),$$
where
$$l_{uv}=\begin{cases}
 \ 1, \ \ \ \ \ \ \ \ \ \ \ \ \ \text{if \ $u=v$};\\
 -1/\sqrt{d_ud_v}, \ \ \text{if $u\thicksim v$};\\
 \ 0, \ \ \ \ \ \ \ \ \ \ \ \ \ \text{otherwise}.
 \end{cases}$$
For brevity, the normalized Laplacian eigenvalues are written as $\mathcal{L}$-eigenvalues. It is well known that the least $\mathcal{L}$-eigenvalue of a connected graph is $0$ with multiplicity $1$ (see \cite{Chung}). Then let the $\mathcal{L}$-eigenvalues of a graph $G$ be
$$\rho_1(G)\geq \rho_2(G)\geq \cdots \geq \rho_{n-1}(G)> \rho_{n}(G)=0.$$

The normalized Laplacian eigenvalues of graphs have been studied intensively (see for example \cite{vanDam6,Guo,Braga,Guo1,Das,Huang2,Sun1,Sun2}), as it reveals not only some structural properties but also some relevant dynamical aspects (such as random walk) of graphs \cite{Chung}. Recently, the multiplicity of the normalized Laplacian eigenvalues attracts much attention. Van Dam and Omidi \cite{vanDam6} determined the graphs with some normalized Laplacian eigenvalue of multiplicity $n-1$ and $n-2$, respectively. Tian $et\ al.$ \cite{Tian,Tian1} characterized some families of graphs of $\mathcal{G}(n, n-3)$, but the graphs with $\rho_{n-1}(G)\neq 1$ and $\nu(G)=diam(G)= 2$ that contain induced $P_4$ are not considered, which is the last remaining case. Sun and Das \cite{Sun3} presented the graphs of $\mathcal{G}(n, n-3)$ with $m(\rho_{n-1}(G))=n-3$ and $m(\rho_{n-2}(G))=n-3$ respectively, and gave the following problem.
\vskip 2 mm
\noindent{\bf Problem\cite{Sun3}:}\ {\sf Is it true that there exists no connected graph with $m(\rho_{1}(G))=n-3\ (n\geq 6)$ and $\nu(G)=2$ ?}
\vskip 2 mm
\noindent
To answer the above problem, it is urgent to complete the characterization of all the graphs in $\mathcal{G}(n, n-3)$. Note that the authors of \cite{Tian,Tian1} have obtained the following results. 

\begin{theorem} \cite{Tian,Tian1}\ Let $G\in \mathcal{G}(n, n-3)$ be a graph of order $n\geq 5$. Then
\end{theorem}
\begin{spacing}{1}
\begin{enumerate}
\item[(i)]  \textit{ $\rho_{n-1}(G)=1$ if and only if $G$ is a complete tripartite graph $K_{a,b,c}$ or $K_n-e$, where $K_n-e$ is the graph obtained from the complete graph $K_n$ by removing an edge.}
\item[(ii)]  \textit{ $\rho_{n-1}(G)\neq 1$ and $\nu(G)\neq 2$ if and only if $G\in \{G_1, G_2, G_3\}$ (see Fig. 1).}
\item[(iii)]  \textit{$\rho_{n-1}(G)\neq 1$, $\nu(G)= 2$ and diam(G)=3 if and only if $G=G_4$ (see Fig. 1).}
\item[(iv)]  \textit{$G$ is a cograph with $\rho_{n-1}(G)\neq 1$ and $\nu(G)= 2$ if and only if $G=G_5$ (see Fig. 1).}
\end{enumerate}
\end{spacing}

Hence, to characterize all graphs of $\mathcal{G}(n, n-3)$ and to address the above problem in \cite{Sun3}, it suffices to consider the graphs that contain induced path $P_4$ with $\rho_{n-1}(G)\neq 1$ and $\nu(G)=diam(G)= 2$. Here, we obtain the following conclusion.

\begin{theorem}\ Let $G\in \mathcal{G}(n, n-3)$ be a graph of order $n\geq 5$. Then $G$ contains induced path $P_4$ with $\rho_{n-1}(G)\neq 1$ and $\nu(G)=diam(G)= 2$ if and only if $G$ is the cycle $C_5$.
\end{theorem}

\begin{remark}\ Combining Theorems 1.1 and 1.2, all graphs of $\mathcal{G}(n, n-3)$ are determined. As a result, the above problem in \cite{Sun3} is answered, that is, there is no connected graph with $m(\rho_{1}(G))=n-3\ (n\geq 6)$ and $\nu(G)=2$. Now, we can also confirm the uncertain result in \cite{Sun3} that if $G$ is the graph with $m(\rho_{1}(G))=n-3$ then $G$ is determined by its normalized Laplacian spectrum.
\end{remark}

\begin{figure}[htbp]
  \centering
  \setlength{\abovecaptionskip}{0cm} 
  \setlength{\belowcaptionskip}{0pt}
  \includegraphics[width=6.5 in]{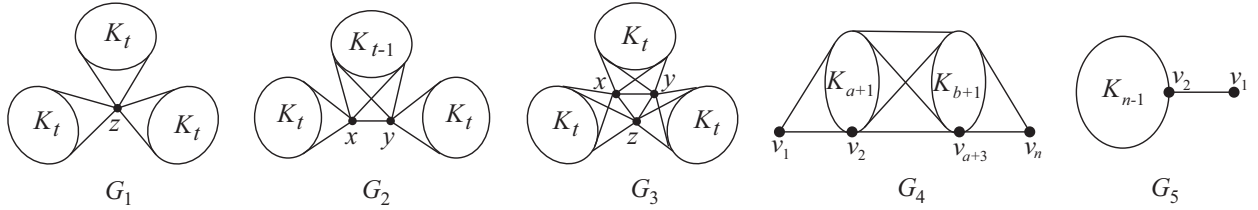}
  \caption{The graphs $G_i\ (1\leq i\leq 5)$.}
\end{figure}

The rest of the paper is arranged as follows. In Section 2, some lemmas are introduced. In Section 3, the proof of Theorem 1.2 is presented.

\section{ Preliminaries }

For brevity, let $\Omega$ be the set of graphs of $\mathcal{G}(n, n-3)$ that contain  induced path $P_4$ with $\rho_{n-1}(G)\neq 1$ and $\nu(G)=diam(G)= 2$. We always let  $\theta$ be the $\mathcal{L}$-eigenvalue of $G$ with multiplicity $n-3$.
\begin{lemma}\label{Sun1lemma}{\rm\cite{Sun1}}\ Let $G$ be a connected graph with order $n\geq 4$. If $G\ncong K_{p,q}, K_a\vee (n-a)K_1$ ($p+q=n, a\geq 2$), then $\rho_{2}(G)\geq \frac{n-1}{n-2}$.\end{lemma}

\begin{lemma}\label{Guolemma}{\rm\cite{Guo}}\ Let $G\ncong K_n$ be a connected graph with order $n\geq 2$. Then $\rho_{n-1}(G)\leq 1$.\end{lemma}

\begin{lemma}\label{cliquelemma}{\rm\cite{Das}}\  Let $G$ be a graph with $n$ vertices. Let $K=\{v_1, \ldots, v_q \}$ be a clique in $G$ such that $N_G(v_i)-K=N_G(v_j)-K$ $(1 \leq i, j\leq q)$, then $1+\frac{1}{d_{v_i}}$ is an $\mathcal{L}$-eigenvalue of $G$ with multiplicity at least $q-1$.\end{lemma}

\begin{lemma}\label{Tian}{\rm\cite{Tian}}\ Let $G\in\mathcal{G}(n, n-3)$ with $\rho_{n-1}(G)\neq 1$, then $\theta\neq 1$.\end{lemma}

\begin{lemma}\label{Tian1}{\rm\cite{Tian1}}\ Let $G\in \Omega$ with an induced path $P_4=v_1v_2v_3v_4$. Then
\begin{equation}\label{e1}
(1-\theta)^4d_{v_1}d_{v_2}d_{v_3}d_{v_4}-(d_{v_1}d_{v_2}+d_{v_3}d_{v_4}+d_{v_1}d_{v_4})(1-\theta)^2+1=0.
\end{equation}
Moreover, if there is a vertex $u_1$ (resp., $u_2$) such that $u_1v_2v_3v_4$ (resp., $v_1u_2v_3v_4$) is also an induced path, then $d_{v_1}=d_{u_1}$ (resp., $d_{v_2}=d_{u_2}$).
\end{lemma}

The following lemma is useful for us to complete the proof of Theorem 1.2.

\begin{lemma}\label{main}\ Let $G\in \Omega$ and $H_i\ (1\leq i\leq 5)$ be the graphs as shown in Fig. 2. Then the following assertions hold. \end{lemma}
\begin{spacing}{1}
\begin{enumerate}
\item[(i)]  \textit{ If $G$ contains $H_1$ as an induced subgraph, then $1-\theta=-\frac{1}{d_{v_1}}=-\frac{1}{d_{v_5}}.$}
\item[(ii)]  \textit{ If $G$ contains $H_2$ as an induced subgraph, then $1-\theta=-\frac{1}{d_{v_2}}=-\frac{1}{d_{v_5}}.$}
\item[(iii)]\textit{If $G$ contains $H_3$ as an induced subgraph, then
$$(1-\theta)^2d_{v_1}d_{v_2}+(1-\theta)d_{v_4}-1=0.$$}
\item[(iv)]\textit{If $G$ contains $H_4$ as an induced subgraph, then
\begin{equation*}
\begin{cases}
(1-\theta)=-\frac{d_{v_3}+d_{v_5}}{d_{v_3}d_{v_2}}=-\frac{d_{v_2}+d_{v_4}}{d_{v_4}d_{v_5}}\\
(1-\theta)^2d_{v_1}(d_{v_3}+d_{v_5})+(1-\theta)d_{v_3}-1=0.
 \end{cases}
\end{equation*}
}
\item[(v)]  \textit{If $G$ contains $H_5$ as an induced subgraph, then $$(1-\theta)=-\frac{d_{v_2}+2d_{v_4}}{d_{v_4}(d_{v_2}+d_{v_5})}=-\frac{d_{v_3}+2d_{v_1}}{d_{v_1}(d_{v_3}+d_{v_5})}.$$
    }
\item[(vi)]\cite{Tian1}  \textit{ If $G$ contains $H_6$ as an induced subgraph, then $1-\theta=-\frac{1}{d_{v_1}}=-\frac{1}{d_{v_4}}.$}
\end{enumerate}
\end{spacing}

\begin{figure}[htbp]
  \centering
  \setlength{\abovecaptionskip}{0cm} 
  \setlength{\belowcaptionskip}{0pt}
  \includegraphics[width=6 in]{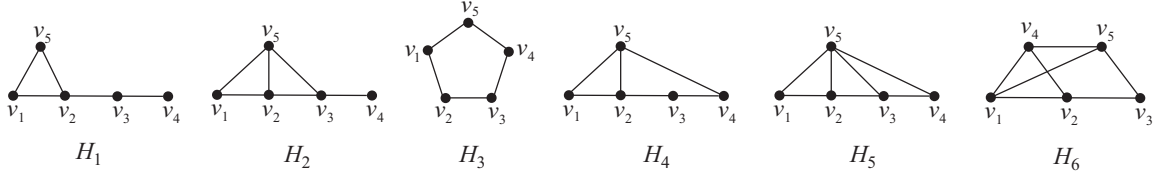}
  \caption{The graphs $H_i\ (1\leq i\leq 6)$.}
\end{figure}

\noindent
{\bf Proof.} \ We first show the proof of assertion (i).

Since $m(\theta)=n-3$, then $r(\mathcal{L}(G)-\theta I)=3$. Lemma \ref{Tian} indicates that $\theta\neq 1$. Denote by $M_1$ the principal submatrix of $\mathcal{L}(G)-\theta I$ indexed by the vertices of $H_1$, then
$$M_1=
\left(
  \begin{array}{ccccc}
    1-\theta & \frac{-1}{\sqrt{d_{v_1}d_{v_2}}} & 0 & 0 & \frac{-1}{\sqrt{d_{v_1}d_{v_5}}}\\
    \frac{-1}{\sqrt{d_{v_1}d_{v_2}}} & 1-\theta & \frac{-1}{\sqrt{d_{v_2}d_{v_3}}} & 0 & \frac{-1}{\sqrt{d_{v_2}d_{v_5}}}\\
    0 & \frac{-1}{\sqrt{d_{v_2}d_{v_3}}} & 1-\theta & \frac{-1}{\sqrt{d_{v_3}d_{v_4}}}& 0 \\
    0 & 0 & \frac{-1}{\sqrt{d_{v_3}d_{v_4}}} & 1-\theta & 0\\
    \frac{-1}{\sqrt{d_{v_1}d_{v_5}}} & \frac{-1}{\sqrt{d_{v_2}d_{v_5}}} & 0 & 0 & 1-\theta\\
  \end{array}
\right).$$
One can easily obtain that the first three rows of $M_1$ are linearly independent (considering the minor indexed by the first three rows and the middle three columns of $M_1$), which yields that the rows $R_{v_1}, R_{v_2}, R_{v_3}$  of $\mathcal{L}(G)-\theta I$ are linearly independent, and then $R_{v_5}$ can be written as a linear combination of $R_{v_1}, R_{v_2}, R_{v_3}$.
Let
\begin{equation}\label{e2}
R_{v_5}=aR_{v_1}+bR_{v_2}+cR_{v_3},
\end{equation}
then
\begin{equation}\label{e3}
\begin{cases}
  a(1-\theta)-\frac{b}{\sqrt{d_{v_1}d_{v_2}}}=\frac{-1}{\sqrt{d_{v_1}d_{v_5}}}\\
\frac{-a}{\sqrt{d_{v_1}d_{v_2}}}+b(1-\theta)-\frac{c}{\sqrt{d_{v_2}d_{v_3}}}=\frac{-1}{\sqrt{d_{v_2}d_{v_5}}}\\
 -\frac{b}{\sqrt{d_{v_2}d_{v_3}}}+c(1-\theta)=0\\
 -\frac{c}{\sqrt{d_{v_3}d_{v_4}}}=0
 \end{cases}
\end{equation}
The fourth equation of (\ref{e3}) implies that $c=0$, and further $b=0$ from the third one. Recalling that $d_{v_1}=d_{v_5}$ by Lemma \ref{Tian1}, then we have $a=1$ from the second of (\ref{e3}), and thus $1-\theta=-\frac{1}{d_{v_1}}=-\frac{1}{d_{v_5}}$ by the first of (\ref{e3}).

For assertion (ii), let
$$M_2=
\left(
  \begin{array}{ccccc}
    1-\theta & \frac{-1}{\sqrt{d_{v_1}d_{v_2}}} & 0 & 0 & \frac{-1}{\sqrt{d_{v_1}d_{v_5}}}\\
    \frac{-1}{\sqrt{d_{v_1}d_{v_2}}} & 1-\theta & \frac{-1}{\sqrt{d_{v_2}d_{v_3}}} & 0 & \frac{-1}{\sqrt{d_{v_2}d_{v_5}}}\\
    0 & \frac{-1}{\sqrt{d_{v_2}d_{v_3}}} & 1-\theta & \frac{-1}{\sqrt{d_{v_3}d_{v_4}}}& \frac{-1}{\sqrt{d_{v_3}d_{v_5}}} \\
    0 & 0 & \frac{-1}{\sqrt{d_{v_3}d_{v_4}}} & 1-\theta & 0\\
    \frac{-1}{\sqrt{d_{v_1}d_{v_5}}} & \frac{-1}{\sqrt{d_{v_2}d_{v_5}}} & \frac{-1}{\sqrt{d_{v_3}d_{v_5}}} & 0 & 1-\theta\\
  \end{array}
\right)$$
be the principal submatrix of $\mathcal{L}(G)-\theta I$ indexed by the vertices of $H_2$.
Similar as above discussion, one can assume that the Eq. (\ref{e2}) still holds. Then
\begin{equation}\label{e4}
\begin{cases}
  a(1-\theta)-\frac{b}{\sqrt{d_{v_1}d_{v_2}}}=\frac{-1}{\sqrt{d_{v_1}d_{v_5}}}\\
\frac{-a}{\sqrt{d_{v_1}d_{v_2}}}+b(1-\theta)-\frac{c}{\sqrt{d_{v_2}d_{v_3}}}=\frac{-1}{\sqrt{d_{v_2}d_{v_5}}}\\
 -\frac{b}{\sqrt{d_{v_2}d_{v_3}}}+c(1-\theta)=\frac{-1}{\sqrt{d_{v_3}d_{v_5}}}\\
 -\frac{c}{\sqrt{d_{v_3}d_{v_4}}}=0.
 \end{cases}
\end{equation}
By the fourth equation of (\ref{e4}), we see $c=0$. Further, recalling that $d_{v_2}=d_{v_5}$ for $H_2$ by Lemma \ref{Tian1}, we get $b=1$ from the third one. Then $a=0$ from the first of (\ref{e4}), and thus $1-\theta=-\frac{1}{d_{v_2}}=-\frac{1}{d_{v_5}}$ by the second one.

For assertion (iii), let
$$M_3=
\left(
  \begin{array}{ccccc}
    1-\theta & \frac{-1}{\sqrt{d_{v_1}d_{v_2}}} & 0 & 0 & \frac{-1}{\sqrt{d_{v_1}d_{v_5}}}\\
    \frac{-1}{\sqrt{d_{v_1}d_{v_2}}} & 1-\theta & \frac{-1}{\sqrt{d_{v_2}d_{v_3}}} & 0 & 0\\
    0 & \frac{-1}{\sqrt{d_{v_2}d_{v_3}}} & 1-\theta & \frac{-1}{\sqrt{d_{v_3}d_{v_4}}}& 0 \\
    0 & 0 & \frac{-1}{\sqrt{d_{v_3}d_{v_4}}} & 1-\theta & \frac{-1}{\sqrt{d_{v_4}d_{v_5}}}\\
    \frac{-1}{\sqrt{d_{v_1}d_{v_5}}} & 0 & 0 & \frac{-1}{\sqrt{d_{v_4}d_{v_5}}} & 1-\theta\\
  \end{array}
\right)$$
be the principal submatrix of $\mathcal{L}(G)-\theta I$ indexed by the vertices of $H_3$.
Clearly, the middle three rows of $M_3$ are linearly independent, which yields that the rows $R_{v_2}, R_{v_3}, R_{v_4}$  of $\mathcal{L}(G)-\theta I$ are linearly independent, and then we set
\begin{equation}\label{e5}
R_{v_1}=aR_{v_2}+bR_{v_3}+cR_{v_4}.
\end{equation}
Applying (\ref{e5}) to the columns of $M_3$, we get
\begin{equation}\label{e6}
\begin{cases}
  \frac{-a}{\sqrt{d_{v_1}d_{v_5}}}=1-\theta\\
a(1-\theta)-\frac{b}{\sqrt{d_{v_2}d_{v_3}}}=\frac{-1}{\sqrt{d_{v_1}d_{v_2}}}\\
 -\frac{b}{\sqrt{d_{v_3}d_{v_4}}}+c(1-\theta)=0\\
 -\frac{c}{\sqrt{d_{v_4}d_{v_5}}}=\frac{-1}{\sqrt{d_{v_1}d_{v_5}}}.
 \end{cases}
\end{equation}
The first and the fourth equations of (\ref{e6}) tell us that $a=-(1-\theta)\sqrt{d_{v_1}d_{v_2}}$ and $c=\sqrt{\frac{d_{v_4}}{d_{v_1}}}$, and further $b=(1-\theta)d_{v_4}\sqrt{\frac{d_{v_3}}{d_{v_1}}}$ from the third one. Taking the values of $a, b, c$ into the second of (\ref{e6}), we derive that $(1-\theta)^2d_{v_1}d_{v_2}+(1-\theta)d_{v_4}-1=0,$ as required.

For assertion (iv), let
$$M_4=
\left(
  \begin{array}{ccccc}
    1-\theta & \frac{-1}{\sqrt{d_{v_1}d_{v_2}}} & 0 & 0 & \frac{-1}{\sqrt{d_{v_1}d_{v_5}}}\\
    \frac{-1}{\sqrt{d_{v_1}d_{v_2}}} & 1-\theta & \frac{-1}{\sqrt{d_{v_2}d_{v_3}}} & 0 & \frac{-1}{\sqrt{d_{v_2}d_{v_5}}}\\
    0 & \frac{-1}{\sqrt{d_{v_2}d_{v_3}}} & 1-\theta & \frac{-1}{\sqrt{d_{v_3}d_{v_4}}}& 0 \\
    0 & 0 & \frac{-1}{\sqrt{d_{v_3}d_{v_4}}} & 1-\theta & \frac{-1}{\sqrt{d_{v_4}d_{v_5}}}\\
    \frac{-1}{\sqrt{d_{v_1}d_{v_5}}} & \frac{-1}{\sqrt{d_{v_2}d_{v_5}}} & 0 & \frac{-1}{\sqrt{d_{v_4}d_{v_5}}} & 1-\theta\\
  \end{array}
\right)$$
be the principal submatrix of $\mathcal{L}(G)-\theta I$ indexed by the vertices of $H_4$.
It is clear that the first three rows of $M_4$ are linearly independent, which indicates that the rows $R_{v_1}, R_{v_2}, R_{v_3}$  of $\mathcal{L}(G)-\theta I$ are linearly independent. Let
\begin{equation}\label{e7}
R_{v_5}=aR_{v_1}+bR_{v_2}+cR_{v_3}.
\end{equation}
Applying (\ref{e7}) to the columns of $M_4$, we have
\begin{equation}\label{e8}
\begin{cases}
  a(1-\theta)-\frac{b}{\sqrt{d_{v_1}d_{v_2}}}=\frac{-1}{\sqrt{d_{v_1}d_{v_5}}}\\
  \frac{-a}{\sqrt{d_{v_1}d_{v_2}}}+b(1-\theta)-\frac{c}{\sqrt{d_{v_2}d_{v_3}}}=\frac{-1}{\sqrt{d_{v_2}d_{v_5}}}\\
\frac{-b}{\sqrt{d_{v_2}d_{v_3}}}+c(1-\theta)=0\\
\frac{-c}{\sqrt{d_{v_3}d_{v_4}}}=\frac{-1}{\sqrt{d_{v_4}d_{v_5}}}\\
\frac{-a}{\sqrt{d_{v_1}d_{v_5}}}-\frac{b}{\sqrt{d_{v_2}d_{v_5}}}=1-\theta.
 \end{cases}
\end{equation}
Combining the last three equations of (\ref{e8}), it follows that
$$a=-\sqrt{d_{v_1}d_{v_5}}(1-\theta)(\frac{d_{v_3}}{d_{v_5}}+1),\ \  b=d_{v_3}(1-\theta)\sqrt{\frac{d_{v_2}}{d_{v_5}}},\ \ c=\sqrt{\frac{d_{v_3}}{d_{v_5}}}.$$
Taking the values of $a, b, c$ into the first and second equations of (\ref{e8}) respectively, one can easily derive that
\begin{equation*}
\begin{cases}
(1-\theta)^2d_{v_1}(d_{v_3}+d_{v_5})+(1-\theta)d_{v_3}-1=0\\
(1-\theta)=-\frac{d_{v_3}+d_{v_5}}{d_{v_3}d_{v_2}}.
 \end{cases}
\end{equation*}
Moreover, by the symmetry between $v_2$ and $v_5$ (resp., $v_3$ and $v_4$) in $H_4$, we can also get $$(1-\theta)=-\frac{d_{v_4}+d_{v_2}}{d_{v_4}d_{v_5}}.$$

At last, we prove assertion (v). Let the principal submatrix of $\mathcal{L}(G)-\theta I$ indexed by the vertices of $H_5$ be $M_5$, then
$$M_5=
\left(
  \begin{array}{ccccc}
    1-\theta & \frac{-1}{\sqrt{d_{v_1}d_{v_2}}} & 0 & 0 & \frac{-1}{\sqrt{d_{v_1}d_{v_5}}}\\
    \frac{-1}{\sqrt{d_{v_1}d_{v_2}}} & 1-\theta & \frac{-1}{\sqrt{d_{v_2}d_{v_3}}} & 0 & \frac{-1}{\sqrt{d_{v_2}d_{v_5}}}\\
    0 & \frac{-1}{\sqrt{d_{v_2}d_{v_3}}} & 1-\theta & \frac{-1}{\sqrt{d_{v_3}d_{v_4}}}& \frac{-1}{\sqrt{d_{v_3}d_{v_5}}} \\
    0 & 0 & \frac{-1}{\sqrt{d_{v_3}d_{v_4}}} & 1-\theta & \frac{-1}{\sqrt{d_{v_4}d_{v_5}}}\\
    \frac{-1}{\sqrt{d_{v_1}d_{v_5}}} & \frac{-1}{\sqrt{d_{v_2}d_{v_5}}} & \frac{-1}{\sqrt{d_{v_3}d_{v_5}}} & \frac{-1}{\sqrt{d_{v_4}d_{v_5}}} & 1-\theta\\
  \end{array}
\right).$$
Similar as the discussion in assertion (iv), the Eq. (\ref{e7}) can still hold. Then applying (\ref{e7}) to the columns of $M_5$, we obtain
\begin{equation}\label{e9}
\begin{cases}
a(1-\theta)-\frac{b}{\sqrt{d_{v_1}d_{v_2}}}=\frac{-1}{\sqrt{d_{v_1}d_{v_5}}}\\
\frac{-b}{\sqrt{d_{v_2}d_{v_3}}}+c(1-\theta)=\frac{-1}{\sqrt{d_{v_3}d_{v_5}}}\\
\frac{-c}{\sqrt{d_{v_3}d_{v_4}}}=\frac{-1}{\sqrt{d_{v_4}d_{v_5}}}\\
\frac{-a}{\sqrt{d_{v_1}d_{v_5}}}-\frac{b}{\sqrt{d_{v_2}d_{v_5}}}-\frac{c}{\sqrt{d_{v_3}d_{v_5}}}=1-\theta.
 \end{cases}
\end{equation}
It follows from the first three equations of (\ref{e9}) that
$$ a=\frac{d_{v_3}}{\sqrt{d_{v_1}d_{v_5}}},\ \ b=\sqrt{d_{v_2}d_{v_3}}((1-\theta)\sqrt{\frac{d_{v_3}}{d_{v_5}}}+\frac{1}{\sqrt{d_{v_3}d_{v_5}}}),\ \ c=\sqrt{\frac{d_{v_3}}{d_{v_5}}}.$$
Taking the values of $a, b, c$ into the last of (\ref{e9}), one can derive
$$(1-\theta)=-\frac{d_{v_3}+2d_{v_1}}{d_{v_1}(d_{v_3}+d_{v_5})}.$$
Furthermore, the symmetry of $H_5$ implies that
$$(1-\theta)=-\frac{d_{v_2}+2d_{v_4}}{d_{v_4}(d_{v_2}+d_{v_5})},$$
as required.

For assertion (vi), one can refer to the process of proving Claim 1 of Lemma 3.2 in \cite{Tian1}. \hfill$\square$

\section{Proof of Theorem 1.2}

\quad   Let $G \in \Omega$, $i.e.,$ $G$ is a graph of $\mathcal{G}(n, n-3)$ containing induced path $P_4$ with $\rho_{n-1}(G)\neq 1$ and $\nu(G)=diam(G)= 2$.
Suppose $m(\theta)=n-3$ in $G$. Now we prove Theorem 1.2.

\vskip 2mm
\noindent
{\bf Proof of Theorem 1.2} \ By direct calculation, the normalized Laplacian spectrum of the cycle $C_5$ is $\{0.691^2, 1.809^2, 0\}$, then it follows that $G \in \Omega$. Thus the sufficiency is clear.

In the following, we present the necessity part.
Suppose that $G \in \Omega$ and $m(\theta)=n-3$ in $G$, then $\theta\neq 1$ from Lemma \ref{Tian}. Denote by $P_4=v_1v_2v_3v_4$ an induced path of $G$. Assume that $U\subseteq V(P_4)$ and $$S_U=\{u\in V(G)\setminus V(P_4): N_G(u)\cap V(P_4)=U\}.$$
It follows from $\nu(G)=2$ that any vertex out of $V(P_4)$ must be adjacent to at least two  of $V(P_4)$ and $S_{\{v_1,v_3\}}=S_{\{v_2,v_4\}}=S_{\{v_2,v_3\}}=\emptyset.$
Further, since $diam(G)=2$, then $d(v_1,v_4)=2$, and thus there exists a vertex, say $v_5$,  adjacent to $v_1$ and $v_4$. Note that $v_5$ maybe belong to $S_{\{v_1,v_4\}}$, $S_{\{v_1,v_2,v_4\}}$, $S_{\{v_1,v_3,v_4\}}$ or $S_{\{v_1,v_2,v_3,v_4\}}$. Accordingly, the remaining proof can be divided into the following cases.

\vskip 2mm
\noindent
{\bf Case 1.} \ Suppose that $v_5\in S_{\{v_1,v_4\}}$, $i.e.,$ $S_{\{v_1,v_4\}}\neq \emptyset$.

We will complete the discussion of this case by the following claims.
\vskip 2mm
\noindent
{\sf Claim 1.1} \ $|S_{\{v_1,v_4\}}|=1$ and $S_{\{v_1,v_2\}}=S_{\{v_3,v_4\}}=\emptyset$.

Suppose that $|S_{\{v_1,v_4\}}|\geq 2$, then all the vertices of $S_{\{v_1,v_4\}}$ induce a clique (otherwise, $\nu(G)\geq 3$, a contradiction). Then $G$ contains an induced subgraph isomorphic to $X_1$ in Fig. 3. Similarly, one can obtain that if $S_{\{v_1,v_2\}}\neq \emptyset$ or $S_{\{v_3,v_4\}}\neq \emptyset$, $G$ also contains an induced subgraph isomorphic to $X_1$.  Since $X_1$ contains $H_1$ as an induced subgraph, then by Lemma \ref{main} (i)
\begin{equation}\label{e10}
1-\theta=-\frac{1}{d_{v_1}}.
\end{equation}
Moreover, $C_5$ ($i.e.,$ $H_3$) is an induced subgraph of $X_1$, then by Lemma \ref{main} (iii)
\begin{equation}\label{e11}
(1-\theta)^2d_{v_1}d_{v_2}+(1-\theta)d_{v_4}-1=0.
\end{equation}
Combining (\ref{e10}) and (\ref{e11}), we get
\begin{equation}\label{e12}
d_{v_2}=d_{v_1}+d_{v_4}.
\end{equation}
It is not hard to see that
\begin{equation}\label{e13}
\begin{cases}
d_{v_2}=|S_{\{v_1,v_2\}}|+|S_{\{v_1,v_2,v_3\}}|+|S_{\{v_2,v_3,v_4\}}|+|S_{\{v_1,v_2,v_4\}}|+|S_{\{v_1,v_2,v_3,v_4\}}|+2\\
d_{v_1}=|S_{\{v_1,v_2\}}|+|S_{\{v_1,v_4\}}|+|S_{\{v_1,v_2,v_3\}}|+|S_{\{v_1,v_2,v_4\}}|+|S_{\{v_1,v_3,v_4\}}|+|S_{\{v_1,v_2,v_3,v_4\}}|+1\\
d_{v_4}=|S_{\{v_3,v_4\}}|+|S_{\{v_1,v_4\}}|+|S_{\{v_2,v_3,v_4\}}|+|S_{\{v_1,v_2,v_4\}}|+|S_{\{v_1,v_3,v_4\}}|+|S_{\{v_1,v_2,v_3,v_4\}}|+1
\end{cases}
\end{equation}
It follows from (\ref{e12}) and (\ref{e13}) that $|S_{\{v_1,v_4\}}|=0$, a contradiction.

\begin{figure}[htbp]
  \centering
  \setlength{\abovecaptionskip}{0cm} 
  \setlength{\belowcaptionskip}{0pt}
  \includegraphics[width=6 in]{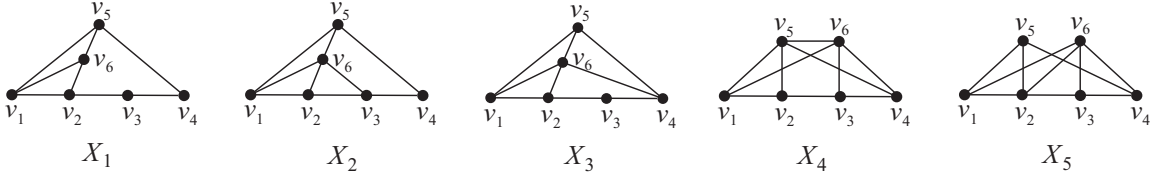}
  \caption{The graphs $X_i\ (1\leq i\leq 5)$.}
\end{figure}

\vskip 2mm
\noindent
{\sf Claim 1.2} \ $S_{\{v_1,v_2,v_3\}}=S_{\{v_2,v_3,v_4\}}=\emptyset$.

It suffices to prove that $S_{\{v_1,v_2,v_3\}}=\emptyset$. Suppose for a contradiction that $S_{\{v_1,v_2,v_3\}}\neq \emptyset$ and $v_6\in S_{\{v_1,v_2,v_3\}}$. If $v_5\nsim v_6$, then the vertices $v_i\ (1\leq i\leq 6)$ induce a subgraph isomorphic to $X_1$. One can obtain a contradiction by similar discussion as above. If $v_5\thicksim v_6$, then $X_2$ in Fig. 3 is an induced subgraph of $G$. Deleting $v_3$ with the incident edges from $X_2$, we also get (\ref{e10}) by Lemma \ref{main} (ii). Analogous discussion as Claim 1.1, the Eq. (\ref{e12}) still holds.
As $|S_{\{v_1,v_4\}}|=1$ and $S_{\{v_1,v_2\}}=S_{\{v_3,v_4\}}=\emptyset$ from Claim 1.1, then
\begin{equation*}
\begin{cases}
d_{v_2}=|S_{\{v_1,v_2,v_3\}}|+|S_{\{v_2,v_3,v_4\}}|+|S_{\{v_1,v_2,v_4\}}|+|S_{\{v_1,v_2,v_3,v_4\}}|+2\\
d_{v_1}=|S_{\{v_1,v_4\}}|+|S_{\{v_1,v_2,v_3\}}|+|S_{\{v_1,v_2,v_4\}}|+|S_{\{v_1,v_3,v_4\}}|+|S_{\{v_1,v_2,v_3,v_4\}}|+1\\
d_{v_4}=|S_{\{v_1,v_4\}}|+|S_{\{v_2,v_3,v_4\}}|+|S_{\{v_1,v_2,v_4\}}|+|S_{\{v_1,v_3,v_4\}}|+|S_{\{v_1,v_2,v_3,v_4\}}|+1,
\end{cases}
\end{equation*}
which, together with (\ref{e12}), yields that $|S_{\{v_1,v_4\}}|=0$, a contradiction.

\vskip 2mm
\noindent
{\sf Claim 1.3} \ $S_{\{v_1,v_2,v_4\}}=S_{\{v_1,v_3,v_4\}}=\emptyset$.

We only need to show $S_{\{v_1,v_2,v_4\}}=\emptyset$. Assume that $S_{\{v_1,v_2,v_4\}}\neq \emptyset$ and $v_6 \in S_{\{v_1,v_2,v_4\}}$, then $v_5\thicksim v_6$ (otherwise $\nu(G)\geq 3$, a contradiction). Thus $X_3$ in Fig. 3 is an induced subgraph of $G$. Removing $v_4$ with the incident edges from $X_3$, we get (\ref{e10}) again by Lemma \ref{main} (ii). By analogous discussion as above, one can easily obtain a contradiction.

\vskip 2mm
\noindent
{\sf Claim 1.4} \ $S_{\{v_1,v_2,v_3, v_4\}}=\emptyset$.

Combining the first three claims, we see that if $S_{\{v_1,v_2,v_3, v_4\}}\neq \emptyset$, then all the vertices of $V(G)\setminus \{v_1,\cdots,v_5\}$ belong to $S_{\{v_1,v_2,v_3, v_4\}}$. Then
\begin{equation}\label{e14}
d_{v_1}=d_{v_2}=d_{v_3}=d_{v_4}.
\end{equation}
Since $G$ now contains $H_5$ (see Fig. 2) as an induced subgraph, then by Lemma \ref{main} (v)
\begin{equation*}
  (1-\theta)=-\frac{d_{v_2}+2d_{v_4}}{d_{v_4}(d_{v_2}+d_{v_5})},
\end{equation*}
which is a rational number. Furthermore, since $G$ contains an induced $C_5$, then by (\ref{e11}) and (\ref{e14}) we derive that
$$(1-\theta)=\frac{-1\pm \sqrt{5}}{2d_{v_1}},$$
which is an irrational number, a contradiction.

The above four claims indicate that if $S_{\{v_1,v_4\}}\neq \emptyset$, then $|S_{\{v_1,v_4\}}|=1$ and $|V(G)|=5$, $i.e.,$ $G$ is the cycle $C_5$.

\vskip 2mm
\noindent
{\bf Case 2.} \ Suppose that $v_5\in S_{\{v_1,v_2,v_4\}}$, $i.e.,$ $S_{\{v_1,v_2,v_4\}}\neq \emptyset$ and $S_{\{v_1,v_4\}}=\emptyset$.

The following claims will help us complete the discussion of this case.
\vskip 2mm
\noindent
{\sf Claim 2.1} \ $S_{\{v_1,v_2\}}=S_{\{v_3,v_4\}}=\emptyset$.

We first demonstrate $S_{\{v_1,v_2\}}=\emptyset$. If $S_{\{v_1,v_2\}}\neq \emptyset$, say $v_6\in S_{\{v_1,v_2\}}$, then by Lemma \ref{main} (i), the Eq. (\ref{e10}) holds. Since $G$ contains $H_4$ (see Fig. 2) as an induced subgraph, then by Lemma \ref{main} (iv)
\begin{equation}\label{e15}
\begin{cases}
(1-\theta)^2d_{v_1}(d_{v_3}+d_{v_5})+(1-\theta)d_{v_3}-1=0\\
(1-\theta)=-\frac{d_{v_2}+d_{v_4}}{d_{v_4}d_{v_5}}.
 \end{cases}
\end{equation}
By (\ref{e10}) and the first equation of (\ref{e15}), we get $d_{v_1}=d_{v_5}$, which implies that $(1-\theta)=-\frac{1}{d_{v_5}}$. Hence, by the second equation of (\ref{e15}), we have $d_{v_2}=0$, a contradiction.

Next, we prove that $S_{\{v_3,v_4\}}=\emptyset$. Suppose that $S_{\{v_3,v_4\}}\neq \emptyset$ and $v_6\in S_{\{v_3,v_4\}}$, then Lemma \ref{main} (i) indicates that
\begin{equation}\label{e16}
  (1-\theta)=-\frac{1}{d_{v_4}}.
\end{equation}
It follows from (\ref{e16}) and the second equation of (\ref{e15}) that $d_{v_2}+d_{v_4}=d_{v_5}$.
Note that any vertex out of $V(P_4)$ must be adjacent to $v_2$ or $v_4$ (thanks to $\nu(G)=2$). Thus
$$d_{v_5}=d_{v_2}+d_{v_4}\geq 3+3+n-6=n,$$
a contradiction. Therefore, $S_{\{v_1,v_2\}}=S_{\{v_3,v_4\}}=\emptyset$, as required.

\vskip 2mm
\noindent
{\sf Claim 2.2} \ $S_{\{v_1,v_2,v_3\}}=S_{\{v_2,v_3,v_4\}}=\emptyset$.

If $S_{\{v_1,v_2,v_3\}}\neq \emptyset$, then $G$ contains $H_2$ and $H_4$ as induced subgraphs. Thus, from Lemma \ref{main} (ii) and (iv),
\begin{equation*}
\begin{cases}
(1-\theta)=-\frac{1}{d_{v_2}}\\
(1-\theta)=-\frac{d_{v_3}+d_{v_5}}{d_{v_3}d_{v_2}},
 \end{cases}
\end{equation*}
which yield that $d_{v_5}=0$, a contradiction.

Similarly, if $S_{\{v_2,v_3,v_4\}}\neq \emptyset$, then from Lemma \ref{main} (ii) and (iv),
\begin{equation*}
\begin{cases}
(1-\theta)=-\frac{1}{d_{v_3}}\\
(1-\theta)=-\frac{d_{v_3}+d_{v_5}}{d_{v_3}d_{v_2}},
 \end{cases}
\end{equation*}
which yield that $d_{v_3}+d_{v_5}=d_{v_2}$. Notice that any vertex distinct with $v_3$ and $v_5$ must be adjacent to $v_3$ or $v_5$ (thanks to $\nu(G)=2$). Therefore,
$$d_{v_2}=d_{v_3}+d_{v_5}\geq 3+3+n-6=n,$$
a contradiction.

\vskip 2mm
\noindent
{\sf Claim 2.3} \ $S_{\{v_1,v_3,v_4\}}=\emptyset$.

Assume that $S_{\{v_1,v_3,v_4\}}\neq \emptyset$ and $v_6\in S_{\{v_1,v_3,v_4\}}$. If $v_5\thicksim v_6$, then $X_4$ (see Fig. 3) is an induced subgraph of $G$. By observation, $X_4$ contains an induced subgraph (by deleting $v_4$ with incident edges) isomorphic to $H_6$ (see Fig. 2), then from Lemma \ref{main} (vi)
\begin{equation}\label{e17}
(1-\theta)=-\frac{1}{d_{v_5}}.
\end{equation}
Combining (\ref{e17}) and the second equation of (\ref{e15}), we obtain that $d_{v_2}=0$, a contradiction.

If $v_5\nsim v_6$, then the principal submatrix, say $M_6$, of $\mathcal{L}(G)-\theta I$ indexed by $\{v_1,\cdots, v_6\}$ can be written as the following block form
$$M_6=
\left(
  \begin{array}{cc}
    M_4 & \alpha \\
    \alpha^T & 1-\theta \\
  \end{array}
\right)
,$$
where $M_4$ has been given in the proof of Lemma \ref{main} (iv) and $$\alpha=(\frac{-1}{\sqrt{d_{v_1}d_{v_6}}}, 0, \frac{-1}{\sqrt{d_{v_3}d_{v_6}}}, \frac{-1}{\sqrt{d_{v_4}d_{v_6}}}, 0)^T,$$
a column vector. Obviously, the Eq. (\ref{e7}) still holds here, and by applying it to the columns of $M_6$, we get the equations of (\ref{e8}) and
\begin{equation}\label{e18}
\frac{-a}{\sqrt{d_{v_1}d_{v_6}}}-\frac{c}{\sqrt{d_{v_3}d_{v_6}}}=0.
\end{equation}
Then from (\ref{e18}) and the values of $a$ and $c$ obtained before
$$a=-\sqrt{d_{v_1}d_{v_5}}(1-\theta)(\frac{d_{v_3}}{d_{v_5}}+1),\ \ c=\sqrt{\frac{d_{v_3}}{d_{v_5}}},$$
it follows that
$\sqrt{\frac{d_{v_1}}{d_{v_5}}}=\sqrt{d_{v_1}d_{v_5}}(1-\theta)(\frac{d_{v_3}}{d_{v_5}}+1)$, which yields $(1-\theta)>0$, contradicting with Lemma \ref{main} (iv).

\vskip 2mm
\noindent
{\sf Claim 2.4} \ $S_{\{v_1,v_2,v_3,v_4\}}=\emptyset$.

Assume that there is a vertex, say $v_6$, in $S_{\{v_1,v_2,v_3,v_4\}}$. If $v_5\nsim v_6$, then $X_5$ (see Fig. 3) is an induced subgraph of $G$. Deleting the vertex $v_3$ with the incident edges from $X_5$, the resultant graph is isomorphic to $H_6$ in Fig. 2. Then from Lemma \ref{main} (vi), $(1-\theta)=-\frac{1}{d_{v_2}}$, which together with
$(1-\theta)=-\frac{d_{v_3}+d_{v_5}}{d_{v_3}d_{v_2}}$ (thanks to Lemma \ref{main} (iv))indicates that $d_{v_5}=0$, a contradiction.

If $v_5\thicksim v_6$, then the principal submatrix, say $M_7$, of $\mathcal{L}(G)-\theta I$ indexed by $\{v_1,\cdots, v_6\}$ can be written as
$$M_7=
\left(
  \begin{array}{cc}
    M_4 & \beta \\
    \beta^T & 1-\theta \\
  \end{array}
\right)
,$$
where $M_4$ is as above and
$$\beta=(\frac{-1}{\sqrt{d_{v_1}d_{v_6}}}, \frac{-1}{\sqrt{d_{v_2}d_{v_6}}}, \frac{-1}{\sqrt{d_{v_3}d_{v_6}}}, \frac{-1}{\sqrt{d_{v_4}d_{v_6}}}, \frac{-1}{\sqrt{d_{v_5}d_{v_6}}})^T,$$
a column vector.
Applying (\ref{e7}) to the last column of $M_7$, we get
\begin{equation}\label{e19}
\frac{a}{\sqrt{d_{v_1}}}+\frac{b}{\sqrt{d_{v_2}}}+\frac{c}{\sqrt{d_{v_3}}}=\frac{1}{\sqrt{d_{v_5}}},
\end{equation}
which together with the second equation of (\ref{e8}) yields that
$$b(\frac{1}{\sqrt{d_{v_2}}}+\sqrt{d_{v_2}}(1-\theta))=0.$$
Since $b\neq 0$ obtained before, then we have
$(1-\theta)=-\frac{1}{d_{v_2}}$, and thus $d_{v_5}=0$  (thanks to $(1-\theta)=-\frac{d_{v_3}+d_{v_5}}{d_{v_3}d_{v_2}}$ in Lemma \ref{main} (iv)), a contradiction.

In this case, combining Claims 2.1-2.4, we see that all the vertices out of $V(P_4)$ belong to $S_{\{v_1,v_2,v_4\}}$. Furthermore, it is obvious that $S_{\{v_1,v_2,v_4\}}$ induces a clique of $G$, as $\nu (G)=2$. Hence, the structure of $G$ is clear now, and $d_{v_2}=n-2$, $d_{v_3}=2$, $d_{v_4}=n-3$ and $d_{v_5}=n-2$. From Lemma \ref{main} (iv),
\begin{equation*}
\begin{array}{rcl}
  (1-\theta)&=&-\frac{d_{v_2}+d_{v_4}}{d_{v_4}d_{v_5}}=\frac{-(2n-5)}{(n-2)(n-3)}\\
  &=& -\frac{d_{v_3}+d_{v_5}}{d_{v_3}d_{v_2}}=\frac{-n}{2(n-2)},
\end{array}
\end{equation*}
which implies that $n=5$, $i.e.,$ $G=H_4$. However, $H_4\notin \Omega$ by direct calculation. Therefore, $S_{\{v_1,v_2,v_4\}}=\emptyset$, and by symmetry we get $S_{\{v_1,v_3,v_4\}}=\emptyset$.

\vskip 2mm
\noindent
{\bf Case 3.} \ Suppose that $v_5\in S_{\{v_1,v_2,v_3,v_4\}}$ ($i.e.,$ $S_{\{v_1,v_2,v_3,v_4\}}\neq \emptyset$) and $S_{\{v_1,v_4\}}=S_{\{v_1,v_2,v_4\}}=S_{\{v_1,v_3,v_4\}}=\emptyset$.

If this is the case, then the vertices of $V(G)\setminus \{V(P_4)\cup S_{\{v_1,v_2,v_3,v_4\}}\}$ belong to $S_{\{v_1,v_2\}}$, $S_{\{v_3,v_4\}}$, $S_{\{v_1,v_2,v_3\}}$ or $S_{\{v_2,v_3,v_4\}}$. Then we have the following claims.

\vskip 2mm
\noindent
{\sf Claim 3.1} \ $S_{\{v_1,v_2\}}=S_{\{v_3,v_4\}}=\emptyset$.

It suffices to show that $S_{\{v_1,v_2\}}=\emptyset$. Suppose $S_{\{v_1,v_2\}}\neq \emptyset$ and $v_6\in S_{\{v_1,v_2\}}$, then
\begin{equation}\label{e20}
(1-\theta)=-\frac{1}{d_{v_1}}
\end{equation}
from Lemma \ref{main} (i). Since
\begin{equation}\label{e21}
(1-\theta)=-\frac{d_{v_2}+2d_{v_4}}{d_{v_4}(d_{v_2}+d_{v_5})}=-\frac{d_{v_3}+2d_{v_1}}{d_{v_1}(d_{v_3}+d_{v_5})}
\end{equation}
from Lemma \ref{main} (v), then by (\ref{e20}) and (\ref{e21}) we derive that
\begin{equation*}
\begin{cases}
d_{v_2}d_{v_1}+2d_{v_1}d_{v_4}=d_{v_2}d_{v_4}+d_{v_4}d_{v_5}\\
2d_{v_1}=d_{v_5},
 \end{cases}
\end{equation*}
which implies that $d_{v_1}=d_{v_4}$. As $G$ contains an induced $P_4$, then the equation (\ref{e1}) holds from Lemma \ref{Tian1}. It follows from (\ref{e1}), (\ref{e20}) and $d_{v_1}=d_{v_4}$ that
\begin{equation}\label{e22}
d_{v_1}(d_{v_2}+d_{v_3})=d_{v_2}d_{v_3}.
\end{equation}
Now we say that $S_{\{v_2,v_3,v_4\}}=\emptyset$, otherwise $(1-\theta)=-\frac{1}{d_{v_3}}$ from Lemma \ref{main} (ii), and then $d_{v_1}=d_{v_3}$ from (\ref{e20}). Thus the equation (\ref{e22}) can be simplified as $d_{v_3}=0$, a contradiction. Analogously, one can derive that $S_{\{v_1,v_2,v_3\}}=\emptyset$. As a result, $d_{v_2}=d_{v_1}+1$. Recalling that $d_{v_1}=d_{v_4}$, then $|S_{\{v_1,v_2\}}|=|S_{\{v_3,v_4\}}|$, and thus $d_{v_2}=d_{v_3}$. Considering (\ref{e22}) again, one can obtain that $d_{v_1}=1$, a contradiction.

\vskip 2mm
\noindent
{\sf Claim 3.2} \ $S_{\{v_1,v_2,v_3\}}=S_{\{v_2,v_3,v_4\}}=\emptyset$.

It suffices to show that $S_{\{v_1,v_2,v_3\}}=\emptyset$. Suppose on the contrary that $S_{\{v_1,v_2,v_3\}}\neq \emptyset$, then $d_{v_2}=d_{v_3}$ by observation. From Lemma \ref{main} (ii), we have $1-\theta=-\frac{1}{d_{v_2}}$. Thus the equation (\ref{e1}) of Lemma \ref{Tian1} can be simplified as
\begin{equation}\label{e23}
d_{v_1}+d_{v_4}=d_{v_2}.
\end{equation}
It is not hard to see that
\begin{equation*}
\begin{cases}
d_{v_1}=|S_{\{v_1,v_2,v_3\}}|+|S_{\{v_1,v_2,v_3,v_4\}}|+1\\
d_{v_4}=|S_{\{v_2,v_3,v_4\}}|+|S_{\{v_1,v_2,v_3,v_4\}}|+1\\
d_{v_2}=|S_{\{v_1,v_2,v_3\}}|+|S_{\{v_2,v_3,v_4\}}|+|S_{\{v_1,v_2,v_3,v_4\}}|+2,
 \end{cases}
\end{equation*}
which together with (\ref{e23}) implies that $|S_{\{v_1,v_2,v_3,v_4\}}|=0$, a contradiction. Therefore, the results of Claim 3.2 hold.

Now we are in a position to complete Case 3. Combining Claims 3.1 and 3.2, we see that all vertices out of $V(P_4)$ belong to $S_{\{v_1,v_2,v_3,v_4\}}$. Then $d_{v_2}=d_{v_1}+1$.
We claim that $S_{\{v_1,v_2,v_3,v_4\}}$ induces a clique of $G$. Otherwise, there exist two vertices, say $v_5$ and $v_6$, of $S_{\{v_1,v_2,v_3,v_4\}}$, which are not adjacent. Then the subgraph induced by $\{v_1,v_2,v_4,v_5,v_6\}$ is isomorphic to $H_6$ in Fig. 2. Hence by Lemma \ref{main} (vi), $1-\theta=-\frac{1}{d_{v_1}}=-\frac{1}{d_{v_2}}$, which indicates that $d_{v_2}=d_{v_1}$, contradicting with $d_{v_2}=d_{v_1}+1$.
As a result, if $|S_{\{v_1,v_2,v_3,v_4\}}|\geq 3$, then $1+\frac{1}{n-1}$ is an $\mathcal{L}$-eigenvalue of $G$ with multiplicity at least 2. Since $G\neq K_n$, then $\rho_{n-1}(G)\leq 1$ by Lemma \ref{Guolemma}. Clearly, $G\ncong K_{p,q}, K_a\vee (n-a)K_1$, then $\rho_{2}(G)\geq \frac{n-1}{n-2}$ by Lemma \ref{Sun1lemma}.
Noting that $\rho_{n}(G)=0$, we obtain $G \notin \Omega$. For the case of $|S_{\{v_1,v_2,v_3,v_4\}}|\leq 2$, one can get $G \notin \Omega$ by direct calculation.

The necessity part can be proved by Cases 1-3, and then the proof is completed. \hfill$\square$

\vskip 3mm
\noindent
{\large\bf Acknowledgements}\\
The authors thank the anonymous referees for their valuable comments of this paper. This work is supported by the Natural Science Foundation of Shandong Province (No. ZR2019BA016).

\end{document}